\documentclass[english,reqno]{amsart}
\usepackage[margin=1in]{geometry}

\usepackage{relsize}
\newcommand{\p}{\partial}

\usepackage{preamble}

\newcommand{\FF}{M_0}
\newcommand{\EE}{M_1}
\newcommand{\GG}{M_2}

\begin{document}

\title{Reconstruction of the Doping Profile in Vlasov-Poisson}

\author{Ru-Yu Lai}
\address{School of Mathematics, University of Minnesota, Minneapolis, MN 55455, USA.}
\email{rylai@umn.edu}
\author{Qin Li} 
\address{Mathematics Department, University of Wisconsin-Madison, 480 Lincoln Dr., Madison, WI 53705 USA.}
\email{qinli@math.wisc.edu}
\author{Weiran Sun}
\address{Department of Mathematics, Simon Fraser University, 8888 University Dr., Burnaby, BC V5A 1S6, Canada}
\email{weirans@sfu.ca}
\date{}
%
%
%
\begin{abstract}
We study the inverse problem of recovering the doping profile in the stationary Vlasov-Poisson equation, given the knowledge of the incoming and outgoing measurements at the boundary of the domain. This problem arises from identifying impurities in the semiconductor manufacturing. Our result states that, under suitable assumptions, the doping profile can be uniquely determined through an asymptotic formula of the electric field that it generates.
\end{abstract}
%
\maketitle
\section{Introduction}
In this paper, we investigate an inverse problem arising from the semiconductor industry. The aim here is to determine the material property (doping profile) of a semiconductor device from external measurements of device characteristics. This problem is known as the \textit{inverse doping problem} and has recently been numerically and analytically studied for various settings. 

In semiconductor manufacturing, characterizing the doping profile of ionized impurity atoms is important for identifying causes of failure and advancing designs. Technologies such as Scanning Capacitance Microscopy (SCM) have been developed to image spatial variations, enabling the mapping of carrier concentration on non-uniformly doped samples. Mathematically, the distributions of electrons and holes in a semiconductor are typically governed by the  drift-diffusion-Poisson system, where the advection term arises from the solution to a Poisson equation which describes the electrostatic potential and uses the doping profile as a source term. In its mathematical formulation, scanning the semiconductor is an inverse problem: One imposes boundary conditions for the  drift-diffusion-Poisson system and makes measurements at the boundary to infer the doping profile in the Poisson equation. Numerous numerical experiments have been conducted~\cite{Kennedy,burger2004inverse,Burger_optimal_design,Hinze_Pinnau,AliTorcicolloVessella,Fang,Leit_o_2006,Cheng_Gamba_Ren}, mostly built upon PDE-constrained optimization approach with an added Tikhonov regularization~\cite{burger2004inverse}.

We note, however, that characterizing the dynamics of holes ($p$) and electrons ($n$) by the drift-diffusion-Poisson system is a relatively crude approximation. These equations are valid only in the macroscopic diffusion limit~\cite{degond_diff,Markowich_book,charge_transport_Jerome}, where the interactions between particles are assumed strong within the domain of interest. As the sizes of the semiconductor devices decrease, the chemical and physical problems enter the nanoscale regimes where the diffusion limit becomes invalid. On this level of fine resolution, the kinetic formulation using the Vlasov-Poisson (VP) system is more accurate~\cite{Carrillo,Poupaud}. Similar to the diffusion limit, the mathematical problem in the inverse problem is to reconstruct the source term in the Poisson equation in this VP system.

Despite the advancement on the technology level, the theoretical understanding of this inverse problem is rather limited except for a few works~\cite{Fang} in a slightly different setting. In particular, the well-posedness result for the VP system in the inverse setting is largely unknown: What collection of data are sufficient to uniquely reconstruct the doping profile? In the current paper, we give a mathematical justification in the kinetic regime. More specifically, we claim that if the boundary source is a thin beam of particles injected in the semiconductor, then by measuring the outgoing location/speed of the particles (either holes or electrons), one can uniquely determine the doping profile.

\subsection{Mathematical problem setup and main result}

As we discussed above, there are several different mathematical models considered in the early literature to describe the transport of electrons in the semiconductor devices, including the drift-diffusion-Poisson system and the Vlasov-Poisson system.
Most of the existing works have contributed to the study of the former system, which is an asymptotic model of the latter one.  
The result for the latter model is relatively sparse.
In the present work, we aim to develop an analytical method to recover the doping profile from the Vlasov-Poisson system with given data.


Let $\Omega \subseteq \R^d$ with $d\geq 2$ 
be the device domain. Suppose that $\Omega$ is bounded and strictly convex with its boundary $\del \Omega$ being $C^2$.   
The semiconductor device we consider is modelled by the Vlasov-Poisson system:
\begin{align} \label{eq:main}
\left\{ \begin{array}{ll}
  v \cdot \nabla_x f(x,v) + E_f(x) \cdot \nabla_x f(x,v) = 0,  & \hbox{ in } \Omega\times \R^d,\\
  -\Delta \phi(x) =  N (x)- \rho_f(x),   & \hbox{ in } \Omega,\\ 
  \phi |_{\del \Omega} (x)= 0 ,\\
    f |_{\Gamma_-}(x,v) = \psi(x, v),\\
\end{array}\right.
\end{align}
where $f$ represents the distribution of electrons in the semiconductor at position $x$ and with velocity $v$. The electric field $E_f(x)$ is generated by the potential $\phi$ and satisfies
\begin{align*}
   E_f(x) = \nabla \phi(x). 
\end{align*}
The source term of the Poisson equation for $\phi$ includes two parts: 
the \textit{doping profile} $N$ given by the semiconductor itself and 
the self-generating part with
$$
    \rho_f(x)=\int_{\R^d} f(x,v)\dv,
$$
where 
$N$ is
the critical characteristic of the system and usually cannot be measured in a nondestructive way.
%
We denote the outgoing and incoming subsets of $\partial\Omega\times\R^d$ by $\Gamma_+$ and $\Gamma_-$ respectively: 
\begin{align*}
\Gamma_{\pm}=\{(x,v)\in\partial\Omega \times\R^d:\,\pm v\cdot n_x>0\},
\end{align*}
where $n_x$ is the unit outward normal vector at $x\in\partial\Omega$.
We impose the zero boundary condition for $\phi$ and the incoming boundary condition $\psi$ for $f$ on the incoming boundary $\Gamma_-$.

For $(x_0, p_0) \in \Gamma_-$, we denote the set 
\begin{align*}
    \mathcal{S}_{x_0,p_0}(\Gamma_-):=\Bigg\{\psi \in W^{1,\infty}(\Gamma_-):\, &0\leq \psi\leq 1,\quad \norm{\nabla_{x, v} \psi}_{L^\infty_{x, v}} \leq
   c_0'|p_0|,\\
   &\quad \Supp(\psi)\subset \vpran{B\left(x_0, {c_0\over|p_0|}\right) \times B\left(p_0, {c_0\over|p_0|}\right)} \cap \Gamma_-\Bigg\},
\end{align*}
where the constants $c_0, c_0'$ depend only on $d,\,\Omega$, and $N$. 
Here $B(x,r)$ is the open ball centered at $x$ with radius $r$.
Based on Theorem~\ref{thm:well-posedness}, 
for a given $\psi\in \mathcal{S}_{x_0,p_0}(\Gamma_-)$ with sufficiently large $|p_0|$, there exists a unique solution $f\in W^{1, \infty}(\Omega\times\R^d)$. Hence we can define the Albedo operator $\Lambda_N: \mathcal{S}_{x_0,p_0}(\Gamma_-)\rightarrow L^{\infty}(\Gamma_+)$ by  
\begin{align*}
\Lambda_N\,:\quad f|_{\Gamma_-} \to f|_{\Gamma_+},
\end{align*}
which maps the incoming data $f|_{\Gamma_-}$ to the outgoing one $f|_{\Gamma_+}$.



Our main result states that the doping profile $N$ can be uniquely determined from the Albedo operator ~$\Lambda_N$.  
 \begin{thm}\label{thm:uniqueness}
    Let $\Omega \subseteq \R^d$ with $d\geq 2$ be an open, bounded, and strictly convex domain with a $C^2$ boundary. Let $\Lambda_{N_1}$ and $\Lambda_{N_2}$ be the Albedo operators associated with $N_1, N_2\in W^{1,\infty}(\Omega)$ in \eqref{eq:main} respectively. 
    Suppose that for any $(x_0, p_0) \in \Gamma_-$ with sufficiently large $|p_0|$, we have
\begin{align*}
    \Lambda_{N_1}(\psi)=\Lambda_{N_2}(\psi)\quad\hbox{ for all }\quad\psi\in \mathcal{S}_{x_0,p_0}(\Gamma_-).
\end{align*}
Then $N_1 = N_2 \hbox{ in }\Omega$. 
\end{thm}

The study of inverse problems for kinetic equations has attracted great attention due to its wide applications, such as medical image, plasma physics, atmospheric science and many others. In particular, the inverse problem for the radiative transfer equation (RTE), also known as the linear Boltzmann equation, is well-studied both analytically and numerically. In the past few decades, there have been many fundamental works contributed to the determination of absorption and scattering coefficients in the RTE from the Albedo operator. One crucial tool is based on the analysis of the singularities of the kernel of the Albedo operator by constructing a special solution concentrated at a given point $(x_0,p_0)\in \Gamma_-$. See for instance \cite{CS1, CS2, CS3, CS98, SU2d} (the uniqueness result) and \cite{Bal14, Bal10, Bal18, lai_inverse_2019, Wang1999, ZhaoZ18} (the stability estimate result) for the demonstration of this singular decomposition method. For more details, we refer to review papers \cite{bal_inverse_2009, Kuireview, Stefanov2003}.

Sharing some similarity to the method above, the main strategy for recovering $N$ in the current paper involves injecting a thin beam of electrons concentrated in both location and velocity into the system.
When the speed $|p_0|$ is large, this type of incoming boundary data allows us to establish the well-posedness for the VP system using the fixed-point theorem. Moreover, it helps to decouple the two equations and transforms the inverse problem for the nonlinear VP system into solving the X-ray transform.

Finally, we note that while the problem is motivated by the study of semiconductors, the mathematical derivation is general. For example, 
it can be applied to recover the inter-particle potential in problems arising from material sciences. See Remark~\ref{rmk:generalization} at the end of Section~\ref{sec:reconstruction}. 



The paper is structured as follows. In Section~\ref{sec:well-posedness}, we introduce several notations and discuss properties of the trajectories of particles. In addition, we apply the fixed point theorem to show the well-posedness of the forward problem. Section~\ref{sec:inverse problems} focuses on how to recover the doping profile from the given data.

\section{The forward problem}\label{sec:well-posedness}
\subsection{Notation}
Throughout the paper, we will use the following notations:
Let $Y$ be one of $\Omega\times\R^d$ and $\Gamma_\pm$. Fix $1\leq p\leq \infty$, for $f\in L^p (Y)$ and $g\in W^{1,p} (Y)$, their norms are simply denoted by 
\begin{align*}
   \|f\|_{L^p_{x,v}} =\|f\|_{L^p (Y)}\quad\hbox{ and }\quad \|g\|_{W^{1,p}_{x,v}} =\|g\|_{W^{1,p} (Y)}.
\end{align*}
 
  Let $F(x)=(F_1(x),\cdots, F_d(x))$ and $G(x)=(G_1(x),\cdots, G_d(x))$ be vector functions in $\R^d$. The vector and matrix notations are defined as follows:
\begin{align} \label{Notation:vec}
\begin{array}{ll}
    \nabla_x F 
    = \left(\partial_{x_i} F_j \right)_{d\times d},\qquad \nabla_x G 
    = \left(\partial_{x_i} G_j \right)_{d\times d},
    \qquad 
    \nabla_x F \cdot G = \left(\sum_{j} \vpran{\del_{x_i} F_j} G_j  \right)_{d1},
    \\   
    G^T \cdot \nabla_x F  = \left(\sum_{i} G_i \partial_{x_i} F_j \right)_{1d}, \qquad \nabla_x G \cdot \nabla_x F 
    = \left(\sum_k \vpran{\del_{x_i} G_k} \vpran{\partial_{x_k} F_j}\right)_{d\times d}.
\end{array}
\end{align}
We will suppress the sub-index in a gradient and use $\nabla F$ whenever no ambiguity arises. 

For $(x,v)\in \overline\Omega\times\R^d$, the characteristics (trajectories) $(X(s), V(s))\equiv (X(s; x, v), V(s; x, v))$ are determined by the Hamilton ODEs:
\begin{align} \label{eq:ODE}
	\frac{{\rm d} X(s;x,v)}{\ds} = V(s; x, v),
	\quad 
	\frac{{\rm d} V(s;x,v)}{\ds} = E (X(s;x,v)),
\end{align}
for $-\infty<s<\infty$ with the initial data $(X(s;x,v), V(s;x,v)) \big|_{s=0} = (x, v)$.
We define the backward exit time $t_-(x,v)\geq 0$ as
\begin{align*}
  t_-(x,v):=\sup(\{0\}\cup \{s> 0:\, X(\tau;x,v)\in \Omega \quad \hbox{for all }\tau\in (-s,0)\}).
\end{align*}
	Similarly, the forward exit time $t_+(x,v)$ is defined by
\begin{align*}
  t_+(x,v):=\sup(\{0\}\cup\{s>0:\, X(\tau;x,v)\in \Omega \quad \hbox{for all }\tau\in (0,s)\}).
\end{align*}
In particular, they satisfy $t_-(x,v)|_{\Gamma_-}=t_+(x,v)|_{\Gamma_+}=0$. 
Moreover, we define the backward and forward exit points and velocities as follows:
\begin{align*}	
   x_- (x,v):= X(-t_-(x,v);x,v),\quad v_- (x,v):= V(-t_-(x,v);x,v),
\\
   x_+ (x,v):= X(t_+(x,v);x,v),\quad v_+ (x,v):= V(t_+(x,v);x,v),
\end{align*}
which satisfy $(x_-(x, v), v_-(x ,v)) \in\Gamma_-$ and $(x_+(x, v), v_+(x ,v)) \in\Gamma_+$.

\subsection{Properties of the characteristics}
We will establish basic results for the trajectories in Lemma~\ref{lemma:gradient trajectory} - Lemma~\ref{lem:non-trapping}. 
The estimates in Lemma~\ref{lemma:gradient trajectory} for the trajectories $(X,V)$ hold by direct computations. Therefore, we omit its proof here.

\begin{lem}\label{lemma:gradient trajectory}
Suppose $(X, V)$ is a solution to~\eqref{eq:ODE} with a given field $E \in W^{1, \infty}(\Omega)$. Then
	\begin{align*}
		\frac{\rm d}{\ds}
		\left(\begin{array}{c}
				\nabla_{x,v} X(s;x,v)\\
				\nabla_{x,v} V(s;x,v)\\
			\end{array}\right)^T
		=  \left(\begin{array}{c}
			\nabla_{x,v} X(s;x,v)\\
			\nabla_{x,v} V(s;x,v)\\
		\end{array}\right)^T \left(\begin{array}{cc}
			O & \nabla E(X(s; x, v))  \\
			I & O \\
		\end{array}\right)
	\end{align*}
and for any $s \in \R$,
\begin{align*}
	\norm{\nabla_{x,v}X(s;x,v)}_{L^\infty(\Omega \times \R^d)}
      + \norm{\nabla_{x,v}V(s;x,v)}_{L^\infty(\Omega \times \R^d)}
  \leq c_d \, e^{c_d(1+ \norm{\nabla E}_{L^\infty(\Omega)
  })|s|}, 
\end{align*}
where $O=O_{d\times d}$ is the $d \times d$ zero matrix, $I=I_{d\times d}$ is the $d \times d$ identity matrix and 
$$
    \left(\begin{array}{c}
				\nabla_{x,v} X(s;x,v)\\
				\nabla_{x,v} V(s;x,v)\\
			\end{array}\right)^T
   :=\left(\begin{array}{cc}
				\nabla_{x} X(s;x,v)& \nabla_{x} V(s;x,v)\\
				\nabla_{v} X(s;x,v) &\nabla_{v} V(s;x,v)\\
			\end{array}\right)_{2d\times 2d}.
$$
Here the constant $c_d$ denotes a constant that only depends on the dimension $d$.  
\end{lem}

In the following, we state and analyze properties of $x_-$, $v_-$, and $t_-$.

\begin{lem}\label{lemma:ids} Let $n_{x_-}$ be the outward normal at $x_- \in \del\Omega$. 
Then we have
\begin{align} \label{eq:x-backward}
  x_-(x,v) = x-vt_- + \int^{-t_-}_0 \int^s_0 E(X(\tau; x, v)) \dtau \ds 
\end{align}
and 
\begin{align} \label{eq:v-backward}
   v_-(x,v) =  v + \int^{-t_-}_0  E(X(\tau; x, v)) \dtau.
\end{align}
Moreover, for $v_-$ satisfying $n_{x_-}\cdot v_-<0$, it holds that	
\begin{align}\label{eq:x-minus-x}
		\nabla_x x_-(x,v) 
	= I - \nabla_x t_- \otimes  v_- + \int^{-t_-}_0\int^s_0 \nabla_x X(\tau) \cdot\nabla E(X(\tau))  \,d\tau ds,
\end{align}
\begin{align}\label{eq:x-minus-v}
  \nabla_v x_-(x,v) 
  = -t_- I - \nabla_v t_- \otimes v_- + \int^{-t_-}_0\int^s_0 \nabla_v X(\tau) \cdot \nabla E(X(\tau)) \,d\tau ds,
\end{align}
\begin{align}\label{DEF:v diff x}
  \nabla_x v_-(x,v) 
= -\nabla_x t_- \otimes E(x_-)
   + \int^{-t_-}_0 \nabla_x X(\tau) \cdot \nabla E(X(\tau)) \,d\tau,
\end{align}
\begin{align}\label{DEF:v diff v}
  \nabla_v v_-(x,v) 
  = I - \nabla_v t_- \otimes E(x_-)
     + \int^{-t_-}_0 \nabla_v X(\tau) \cdot \nabla E(X(\tau)) \,d\tau,
\end{align}
\begin{align} \label{eq:t-minus-dev-x}
   \nabla_x t_-(x,v) 
   = {1\over |n_{x_-}\cdot v_-|}
      \left[ -n_{x_-} - \int^{-t_-}_0\int^s_0 \nabla_x X(\tau) \cdot \vpran{\nabla E(X(\tau)) \cdot n_{x_-}}  \,d\tau ds\right],
\end{align}
\begin{align} \label{eq:t-minus-dev-v}
	\nabla_vt_-(x,v) = {1\over |n_{x_-}\cdot v_-|}\left[ n_{x_-}t_- - \int^{-t_-}_0\int^s_0 \nabla_v X(\tau) \cdot \vpran{\nabla E(X(\tau)) \cdot n_{x_-}} \,d\tau ds\right],
\end{align}
where we simply denote $t_-=t_-(x,v)$, $x_-=x_-(x,v)$ and $v_-=v_-(x,v)$ for $(x,v)\in \Omega\times\R^d$.
\end{lem}
\begin{proof}
	Recall that $x_-=X(-t_-(x,v);x,v)$, $v_-=V(-t_-(x,v);x,v)$ and $x=X(0;x,v)$, $v=V(0;x,v)$. Equations~\eqref{eq:x-backward}-\eqref{eq:v-backward} follow from a direct integration of~\eqref{eq:ODE} from time $0$ to $-t_-$. To derive the equations \eqref{eq:x-minus-x} and \eqref{DEF:v diff x}, we differentiate \eqref{eq:x-backward} and \eqref{eq:v-backward} with respect to  $x$. Similarly, \eqref{eq:x-minus-v} and \eqref{DEF:v diff v} follow by differentiating \eqref{eq:x-backward} and \eqref{eq:v-backward} with respect to $v$.

	To derive equations~\eqref{eq:t-minus-dev-x}-\eqref{eq:t-minus-dev-v}, let $\xi(x)$ be the $C^2$
function which defines the boundary such that 
\begin{align*}
	\partial\Omega = \{\xi(x)=0\},\quad \Omega = \{\xi(x)>0\}.
\end{align*}
Then $n_x = -\nabla \xi(x)/|\nabla\xi(x)|$ is the unit outward normal at each point $x \in \partial\Omega$. Since $x_-(x,v)\in \partial\Omega$, we have $\xi(x_-(x, v)) = 0$ so that the derivatives of $\xi(x_-(x,v))$ with respect to both variables $x$ and $v$ vanish, that is,
	\begin{align*}
		\nabla_x (\xi(x_-(x, v))) = 0 \quad\hbox{and}\quad \nabla_v (\xi(x_-(x, v)))= 0,
	\end{align*}
	which lead to
	\begin{align*}
	    \nabla_x x_-\cdot 	n_{x_-} = 0\quad\hbox{and}\quad   \nabla_v x_- \cdot n_{x_-} = 0.
	\end{align*}
	Hence multiplying ~\eqref{eq:x-minus-x} and~\eqref{eq:x-minus-v} by $n_{x_-}$, respectively, we obtain ~\eqref{eq:t-minus-dev-x} and~\eqref{eq:t-minus-dev-v}. 
\end{proof}

For simplicity, from now on we will denote
\begin{align} \label{def:e-0}
  \FF := \norm{E}_{L^\infty(\Omega)},
\qquad \hbox{and}\qquad
   \EE := \norm{\nabla E}_{L^\infty(\Omega)}.
\end{align}
Using the formulas in Lemma~\ref{lemma:ids}, we can bound the derivatives of $(t_-, x_-, v_-)$ with respect to $(x, v)$. They will be utilized to give an upper bound of the solution of the VP system in the proof of the well-posedness theorem later.
\begin{lem}\label{lemma:Estimates} 
Suppose that there exist $p_0 \in \R^d$ with $|p_0| \gg 1$ and constants $\alpha_1,\,\alpha_2,\,\beta_0 > 0$ and $0<\delta<1$ such that for all $(x, v)\in\Omega\times\R^d$, 
\begin{align} \label{bound:n-dot-v}
   \alpha_1 |p_0|\leq |v_-| \leq \alpha_2 |p_0|,
\qquad
   0 \leq t_- \leq \frac{\beta_0}{|p_0|} \leq 1 , 
\qquad
    n_{x_-}\cdot \frac{v_-}{|v_-|} < - \delta < 0. 
\end{align}
Then 
\begin{align}\label{EST:diff x}
  & \norm{\nabla_x x_-}_{L^\infty_{x, v}} 
\leq C \left( 1+{1\over \delta}  + {M_1\over \delta|p_0|^2}e^{C(1+M_1)\over |p_0|}\right) ,
\end{align}
\begin{align}\label{EST:diff x2}
  &  \norm{\nabla_v x_-}_{L^\infty_{x, v}} 
\leq  
  C \left( {1\over \delta|p_0|} 
  + \frac{M_1}{\delta |p_0|^2} e^{C(1+M_1)\over |p_0|} \right),
\end{align}
\begin{align}\label{EST:v diff x}
  \norm{\nabla_x v_-}_{L^\infty_{x, v}}  
\leq C\left( {M_0\over \delta|p_0|} + {M_0 M_1\over \delta|p_0|^3} e^{C(1+M_1)\over |p_0|} + { M_1\over |p_0|} e^{C(1+M_1)\over |p_0|}\right),
\end{align}
\begin{align}\label{EST:v diff x2}
  \norm{\nabla_v v_-}_{L^\infty_{x, v}}  
\leq C\left(1+ {M_0\over \delta|p_0|^2} + {M_0 M_1\over \delta|p_0|^3}  e^{C(1+M_1)\over |p_0|} + { M_1\over |p_0|} e^{C(1+M_1)\over |p_0|}\right),
\end{align}
\begin{align}\label{EST:t diff x}
  \norm{\nabla_x t_-}_{L^\infty_{x, v}} 
\leq {C\over  \delta|p_0|} \left(1+ {M_1\over |p_0|^2}  e^{C(1+M_1)\over |p_0|} \right),    
\qquad 
  \norm{\nabla_v t_-}_{L^\infty_{x, v}} 
\leq {C\over  \delta|p_0|^2} \left(1+ {M_1\over |p_0|} e^{C(1+M_1)\over |p_0|} \right),
\end{align}
where the constant $C>0$ depends only on $d,\, \alpha_1,\,\alpha_2$ and $\beta_0$.  
\end{lem}
\begin{proof}
In this proof, the constant $C$ may vary from line to line. Applying Lemma~\ref{lemma:gradient trajectory} to equations~\eqref{eq:t-minus-dev-x} and~\eqref{eq:t-minus-dev-v} with $\EE = \norm{\nabla E}_{L^\infty (\Omega)}$, we bound $\nabla_x t_-$ and $\nabla_v t_-$ as follows:
\begin{align} \label{bd:dev-x-t-minus}
  |\nabla_x t_-(x,v)| 
&\leq
  \frac{c_d}{\alpha_1  \delta|p_0|}
  \vpran{1 + \EE \, c_d \int_{-t_-}^0 \int_s^0  e^{c_d (1+\EE) |\tau|} \dtau \ds}  \nn
\\
&\leq
  \frac{c_d}{\alpha_1 \delta |p_0|}
  \vpran{1 + {t_-^2 \over 2} \EE \, c_d \, e^{c_d (1+\EE) t_-}}  \notag\\
&\leq
  \frac{c_d}{\alpha_1  \delta|p_0|}
  \vpran{1 + \frac{\beta_0^2}{2|p_0|^2} \EE c_de^{\frac{c_d \beta_0}{|p_0|} (1+\EE)}} \notag\\
&\leq  {C\over  \delta|p_0|} \left(1+ {M_1\over |p_0|^2}  e^{C(1+M_1)\over |p_0|} \right),
\end{align}
where the positive constant $C$ depends on $d,\,\alpha_1$ and $\beta_0$. Similarly, we derive that
\begin{align} \label{bd:dev-v-t-minus}
  |\nabla_v t_-(x,v)| 
&\leq
 \frac{c_d }{\alpha_1  \delta|p_0|}
 \vpran{t_- + \frac{\beta_0^2}{2|p_0|^2} \EE c_d\, e^{\frac{c_d\beta_0}{|p_0|} (1+\EE)}}  \notag\\
&\leq
 \frac{c_d  \beta_0}{\alpha_1  \delta|p_0|^2}
 \vpran{1 + \frac{\beta_0}{2|p_0|} \EE c_d\, e^{\frac{c_d\beta_0}{|p_0|} (1+\EE)}} \notag\\
&\leq {C\over  \delta|p_0|^2} \left(1+ {M_1\over |p_0|} e^{C(1+M_1)\over |p_0|} \right).
\end{align}
Thus~\eqref{EST:t diff x} holds. We can then apply~\eqref{EST:t diff x} to~\eqref{eq:x-minus-x}-\eqref{DEF:v diff v} and $\delta < 1$ to deduce the following estimates:
\begin{align*}
  \abs{\nabla_x x_-(x,v)} 
&\leq
  c_d \vpran{1 + \alpha_2 \abs{p_0} \abs{\nabla_x t_-} 
		+ \int_{-t_-}^0 \int_s^0 \abs{\nabla_x X} \abs{\nabla E} \dtau \ds }
\\
& \leq
    c_d\left( 1+ \alpha_2|p_0| \left({C\over  \delta|p_0|} \left(1+ {M_1\over |p_0|^2} e^{C(1+M_1)\over |p_0|} \right)\right) + \frac{\beta_0^2}{2|p_0|^2} \EE c_d \, e^{\frac{c_d \beta_0}{|p_0|} (1+\EE)} \right) 
\\
& =
    c_d\left( 1+ {\alpha_2  C  \over  \delta} \left(1+ {M_1\over |p_0|^2} e^{C(1+M_1)\over |p_0|} \right) + \frac{\beta_0^2}{2|p_0|^2} \EE c_d \, e^{\frac{c_d \beta_0}{|p_0|} (1+\EE)} \right) 
 \\
 &\leq  
    C \left( 1+{1\over \delta}  + {M_1\over \delta|p_0|^2}e^{C(1+M_1)\over |p_0|}  +  {M_1\over |p_0|^2}e^{C(1+M_1)\over |p_0|}\right)
 \\
  &\leq  
    C \left( 1+{1\over \delta}  + {M_1\over \delta |p_0|^2}e^{C(1+M_1)\over |p_0|}\right),
\end{align*}
where the positive constant $C$ depends on $d,\,\alpha_1,\,\alpha_2$, and $\beta_0$.
Using similar arguments, we have
\begin{align*}
  \abs{\nabla_v x_-(x,v)} 
&\leq
  c_d \vpran{t_- + |v_-| \abs{\nabla_v t_-}  
		+ \int_{-t_-}^0 \int_s^0 \abs{\nabla_v X} \abs{\nabla E} \dtau \ds }
\\
& \leq
    C \left( {1\over |p_0|} + {1\over \delta|p_0|} + {M_1\over \delta|p_0|^2}e^{C(1+M_1)\over |p_0|}  + {M_1\over |p_0|^2}e^{C(1+M_1)\over |p_0|}  \right)
\\
& \leq
    C \left( {1\over |p_0|} + {1\over \delta|p_0|} + \frac{M_1}{\delta |p_0|^2} e^{C(1+M_1)\over |p_0|} \right).
\end{align*}
For $\nabla_x v_-$ and $\nabla_v v_-$, we then derive
\begin{align*}
    \abs{\nabla_x v_-(x,v)} 
   &\leq c_d \vpran{\abs{\nabla_x t_-} |E| + \int_{-t_-}^0 \abs{\nabla_x X} \abs{\nabla E} \ds }\\
   &\leq C\left( {M_0\over \delta|p_0|} + {M_0 M_1\over \delta|p_0|^3} e^{C(1+M_1)\over |p_0|} + { M_1\over |p_0|} e^{C(1+M_1)\over |p_0|}\right)
\end{align*}
and 
\begin{align*}
	\abs{\nabla_v v_-(x,v)} 
	&\leq
	c_d \vpran{1 + |\nabla_v t_-| |E|
		+ \int_{-t_-}^0 \abs{\nabla_v X} \abs{\nabla E} \ds}\\
    &\leq C\left(1+ {M_0\over \delta|p_0|^2} + {M_0 M_1\over \delta|p_0|^3}  e^{C(1+M_1)\over |p_0|} + { M_1\over |p_0|} e^{C(1+M_1)\over |p_0|}\right),
\end{align*}
by which we finish the proof.
\end{proof}

Our next step is to show that for the incoming data that we use for the reconstruction, the generated characteristics will not be trapped in the domain $\Omega$, that is, there exists a finite time that the trajectory will reach $\del\Omega$.

We denote the diameter of $\Omega$ by $R=\text{diam}(\Omega)$ and recall that $B(p_0,r)$ represents an open ball centered at $p_0$ with radius $r>0$.
\begin{lem} \label{lem:non-trapping}
Let $E \in W^{1,\infty}(\Omega)$. 
Then the characteristic ODE system:
\begin{align} \label{eq:charact}
\begin{split}
    &\dot X = V(t), 
\quad
   \dot V = E(X(t)),
\\
   (X(0), &V(0)) = (x, v) \in (\del\Omega \times B(p_0, |p_0|/2)) \cap \Gamma_-
\end{split}
\end{align}
has a unique non-trapping solution if $|p_0|$ is large enough (see its explicit bounds in~\eqref{bound:t-plus}). 
\end{lem}
\begin{proof}
The integral form of~\eqref{eq:charact} is
\begin{align} \label{eq:integ-charact}
   X(t) = x + v t + \int_0^t \int_0^s E(X(\tau)) \dtau \ds, 
\qquad
  t \in (0, t_+(x, v)), 
\end{align}
where $t_+(x, v) \in [0, \infty]$ is the forward exit time of the trajectory starting from $(x, v) \in (\del\Omega \times B(p_0, |p_0|/2)) \cap \Gamma_-$. 
Since $E$ is globally Lipschitz on $\Omega$, equation~\eqref{eq:integ-charact} has a unique solution which exists as long as $X$ stays in $\Omega$. We show that $t_+ < \infty$ if $|p_0|$ is large enough, that is, $X$ will leave $\Omega$ in a finite time. Recall that $\FF = \norm{E}_{L^\infty(\Omega)}$. Then by~\eqref{eq:integ-charact},
\begin{align} \label{ineq:X-lower-bd}
   \abs{X (t) - x}
\geq
   |v| t - \frac{\FF}{2} t^2. 
\end{align}
Direct computation gives
\begin{align*}
   \max \left\{|v| t - \frac{\FF}{2} t^2 \right\}
= \frac{|v|^2}{2 \FF}. 
\end{align*}
Therefore, 
if we set  
\begin{align*}
   \frac{|v|^2}{2 \FF}  > R, 
\quad \text{or equivalently,} \quad
  |v| > \sqrt{2 \FF R }, 
\end{align*}
then the right-hand side of~\eqref{ineq:X-lower-bd} will reach $R$ at some $t_0 \in (0, \infty)$, which forces $X$ to exit $\Omega$ before $t = t_0$. For later use, denote $t_R$ as the smaller solution to 
\begin{align*}
	|v| t_R - \frac{M_0}{2} t_R^2 = R. 
\end{align*}
That is, $t_R =  \frac{|v| - \sqrt{|v|^2 - 2 M_0 R}}{M_0} = \frac{2R}{|v| + \sqrt{|v|^2 - 2 M_0R}}$. Since $v\in B(p_0,|p_0|/2)$ implies $|v|>|p_0|/2$, this gives the upper bound of the exit time
\begin{align} \label{bound:t-plus}
	t_+ \leq t_R= \frac{2 R}{|v| + \sqrt{|v|^2 - 2 M_0R}} 
\leq
	\frac{2R}{|v|}
<
	\frac{4R}{|p_0|}, 
\qquad
	\text{if $|p_0| \geq 1 + 2 \sqrt{2M_0R}$.}     
\end{align}
Hence the upper bound of the exit time only depends on $\Omega$ and $|p_0|$ and it decreases when $|p_0| \to \infty$. 
\end{proof}

\subsection{Proof of the well-posedness}
In this section, we will discuss the forward problem for the Vlasov-Poisson system \eqref{eq:main}. The proof is based on the application of the fixed-point argument, which enables us to show that there exists a unique solution for the system under suitable assumptions. 

Recall that $\phi$ is the solution to the Dirichlet problem: $-\Delta\phi=N-\rho_g$ in $\Omega$, $\phi=0$ on $\partial\Omega$. Based on Lemma~4.1, 4.2 and Theorem~4.3 in \cite{Gilbarg}, the problem is uniquely solvable if $N,\, g\in W^{1,\infty}(\Omega)$ and $\Supp(g (x,\cdot))\subseteq  B(p_0, 2\Eps)$ for $\epsilon>0$.
Let $G=G(x,y)$ be the Green's function of the Laplace operator over the domain $\Omega$. Thus we have
$$
    \phi(x)= \int_\Omega G(x,y)(N-\rho_g)(y)\,dy=:G \ast (N-\rho_g)\quad\hbox{for $x\in\Omega$}.
$$
Since the electric field $E_g$ is coupled to the distribution $g$ by the Poisson equation, it satisfies
\begin{align}\label{Def:E}
    E_g(x) = \nabla \phi(x) = \nabla_x\int_\Omega  G(x,y)(N - \rho_g)(y)\,dy=\nabla_x G \ast (N-\rho_g)(x)\quad\hbox{for $x\in\Omega$}.
\end{align}

The purpose of the following two lemmas is to derive upper bounds for the external force $E_g$. 
\begin{lem} \label{lem:Green}
Suppose $\Omega$ is strictly convex with $C^1$ boundary. Then
\begin{itemize}
\item[(a)] If $d \geq 3$, then 
\begin{align}\label{EST:G-d}
	|\nabla_xG(x,y)| \leq {C\over |x-y|^{d-1}},\quad |\nabla_x^2G(x,y)| \leq {C\over |x-y|^{d}},
\end{align} 
where $C$ depends only on $\Omega$ and $d$. 

\item[(b)] If $d = 2$, then 
\begin{align}\label{EST:G-2}
   |\nabla_xG(x,y)| 
\leq 
  {C\over |x-y|} \vpran{1 + \abs{\ln |x - y|}},
\quad
   |\nabla_x^2G(x,y)| 
 \leq {C\over |x-y|^2} \vpran{1 + \abs{\ln |x - y|}},
\end{align} 
where $C$ depends only on $\Omega$. 
\end{itemize}
\end{lem}
The proof of Lemma~\ref{lem:Green} in the case of $d \geq 3$ is shown in \cite{Widman1967}. An argument along a similar line gives the estimates for $d = 2$.


 
As we will see, the bounds on the electric field $E_g$ follow from the bounds of the Green's function. 


\begin{lem}\label{lemma:E}
Let $\epsilon>0$ and $p_0$ be a fixed vector in $\R^d$. Let $N\in W^{1,\infty}(\Omega)$. For $g\in W^{1,\infty}(\Omega \times\R^d)$ satisfying $\Supp(g (x,\cdot))\subseteq  B(p_0, 2\Eps)$, we have
	\begin{align}\label{EST:E}
		\|E_g\|_{L^\infty(\Omega)}\leq 	C \vpran{\norm{N}_{L^\infty(\Omega)} + \Eps^d  \norm{g}_{L^{\infty}_{x,v}} } 
	\end{align}
and 
    \begin{align}\label{EST:gradient E}
	    \|\nabla E_g\|_{L^\infty(\Omega)}
    	\leq C \left(\|N\|_{W^{1,\infty}(\Omega)} + \epsilon^d \norm{g}_{W^{1,\infty} _{x,v}}\right),
    \end{align}
where the constant $C>0$ depends on $\Omega$ and $d$.  
\end{lem}
\begin{proof}
    For all $x\in\Omega$, by Lemma~\ref{lem:Green}, we have
	\begin{align*}
		|E_g(x)|
		&\leq  \int_\Omega |\nabla_xG(x,y)||(N - \rho_g)(y)|\,dy\\
		&\leq  C \|N - \rho_g\|_{L^\infty(\Omega)}\int_\Omega |x-y|^{1-d} \vpran{1 + \abs{\ln |x-y|}} \,dy  \\
		&\leq  C (\|N\|_{L^\infty(\Omega)} + \|\rho_g\|_{L^\infty(\Omega)}),\
	\end{align*}
    where $C$ depends on $\Omega$. Since $g$ is supported in $v$, we have
	\begin{align*}
		\norm{\rho_g}_{L^\infty(\Omega)}
		= \norm{\int_{\R^d} g(x, v) \dv}_{L^\infty(\Omega)}
		= \norm{\int_{B(p_0, 2\Eps)}g(x, v) \dv}_{L^\infty(\Omega)}
		\leq
		c_d\Eps^d \norm{g}_{L^\infty_{x, v}}, 
	\end{align*}
	where $c_d$ only depends on the dimension $d$. Altogether, we have
	\begin{align*}
		\norm{E_g}_{L^\infty(\Omega)}
		\leq
		C  
		\vpran{\norm{N}_{L^\infty(\Omega)} + \Eps^d  \norm{g}_{L^\infty_{x, v}}}. 
	\end{align*}
	Moreover, we also bound $\norm{\nabla E_g}_{L^\infty(\Omega)}$ by
	\begin{align*}
		|\nabla E_g(x)|
		&= \left| \nabla_x^2\int_\Omega  G(x,y) (N - \rho_g)(y) \,dy\right|\\
		&\leq  \left|\int_\Omega  \nabla_x^2 G(x,y) ((N - \rho_g)(y)-(N - \rho_g)(x)) \,dy\right|+ |(N - \rho_g)(x)|\left|\nabla_x^2 \int_\Omega   G(x,y)   \,dy\right|\\
		&=: J_1+J_2.
	\end{align*}
To estimate $J_1$, by Lemma~\ref{lem:Green} again, we can derive that
\begin{align*}
   &\quad \left|\int_\Omega  \nabla_x^2 G(x,y) ((N - \rho_g)(y)-(N - \rho_g)(x)) \,dy\right|
\\
&\leq 
   \int_{\{|x-y|\leq R\}} {C\over|x-y|^{d-1}} {|(N - \rho_g)(y)-(N - \rho_g)(x)| \over |x-y|} \vpran{1 + \abs{\ln |x-y|}} \,dy 
\\
&\leq 
   \|\nabla(N-\rho_g)\|_{L^\infty(\Omega)}\int_{\{|x-y|\leq R\}} {C\over|x-y|^{d-1}}\vpran{1 + \abs{\ln |x-y|}} \,dy 
   \\
& \leq 
  C \vpran{\norm{N}_{W^{1, \infty} (\Omega)} + \norm{\rho_g}_{W^{1, \infty} (\Omega)}}
\leq
  C \vpran{\norm{N}_{W^{1, \infty} (\Omega)} + \Eps^d \norm{g}_{W^{1,\infty}_{x,v}}}, 
\end{align*} 
where $R$ is the diameter of $\Omega$. 
Now for $J_2$, we first let
$$
u(x) := \int_\Omega G(x,y) \,dy.
$$
Then $u$ is a solution to the Laplace equation $-\Delta u=1$ in $\Omega$ with $u=0$ on $\partial\Omega$. This implies that $u\in C^2(\overline\Omega)$ if $\partial\Omega\in C^2$ and has an upper bound
\begin{align*}
    \|u\|_{C^2(\overline\Omega)}\leq C,
\end{align*}
where the constant $C>0$ depends on $\Omega$ \cite{Evans}. Therefore, we derive that 
\begin{align*}
   J_2
\leq 
   C\|N - \rho_g\|_{L^\infty(\Omega)}
\leq
    C \vpran{\norm{N}_{W^{1, \infty} (\Omega)} + \Eps^d \norm{g}_{W^{1,\infty}_{x,v}}}.
\end{align*}
Combining the above estimates for $J_1$ and $J_2$ gives the desired bound for $\nabla E_g$.
\end{proof}

Using the bounds derived above, we now state and prove the well-posedness theorem for system~\eqref{eq:main}. 
To begin, we first denote the set 
$$
   \mathcal{M}:=\{h\in  W^{1,\infty}(\Omega):\, \|h\|_{W^{1,\infty}(\Omega)}\leq m_0\}
$$
for a constant $m_0>0$.
\begin{thm} \label{thm:well-posedness}
Let $\Omega$ be an open bounded and strictly convex domain with $C^2$ boundary. 
Let $(x_0, p_0) \in \Gamma_-$ and parameter $0<\delta_0<1$ such that
\begin{align}\label{delta 0}
    n_{x_0} \cdot \frac{p_0}{|p_0|} = -\delta_0 < 0.
\end{align}
Suppose $N \in \mathcal{M}$. Then there exists 
a sufficiently large $|p_0|>1$ such that for $\psi \in W^{1,\infty}(\partial\Omega\times\R^d)$, which is compactly supported on $\vpran{B(x_0, {c_0\over|p_0|}) \times B(p_0, {c_0\over|p_0|})} \cap \Gamma_-$ and satisfies 
\begin{align} \label{bound:deriv-psi}
   0\leq  \psi \leq 1, \quad \norm{\nabla_{x, v} \psi}_{L^\infty_{x, v}} \leq
   c_0'|p_0| 
\end{align}
with fixed positive constants $c_0$ and $c_0'$, the system \eqref{eq:main} has a unique solution $f\in W^{1,\infty}(\Omega\times\R^d)$ satisfying
$$
    \norm{f}_{L^\infty_{x,v}} \leq \norm{\psi}_{L^\infty_{x, v}}.
$$  
\end{thm}

Note that based on the proof below, the constants $c_0$ and $c_0'$ depend only on $\Omega, \, d$, and $N$.
\begin{proof}
The proof follows from the fixed-point argument. 
We start with specifying a set $\Sigma \subseteq \Omega$ which will be a superset of the support of $f$ in $x$. 
Let $\Sigma_0$ be the collection of points on the lines initiating from $\Supp(\psi)$:
\begin{align*}
   \Sigma_0 
:= \{x \in \bar\Omega \, \big| \,
      x = y + sv \,\, \text{for some $(y,v) \in \Supp(\psi)$ and $s \geq 0$}\}.
\end{align*}
We define the set $\Sigma$ by
\begin{align} \label{def:Sigma}
  \Sigma_\epsilon
:= \left\{x \in \bar\Omega \, \Big| \,\, \text{dist}(x, \Sigma_0) \leq 2\epsilon  \right\}. 
\end{align}
In addition, we also define the solution set $\CalX$ as
\begin{align*}
   \CalX := \left\{f \in  W^{1,\infty}(\Omega \times \R^d) \big| \, \norm{f}_{L^\infty_{x,v}} \leq \norm{\psi}_{L^\infty_{x, v}}, \,\, \norm{\nabla_{x, v} f}_{L^\infty_{x,v}}  \leq \GG, \,\, \Supp (f) \subseteq \Sigma_\epsilon \times B(p_0, 2\epsilon) \right\}, 
\end{align*}
with the constants $\GG$ and $\epsilon$ to be determined later (see~\eqref{def:M-2} and~\eqref{def:epsilon-0}).

We will prove the Lipschitz continuity of the operator $\mathcal{T}:g\mapsto f$ in the $L^\infty(\Omega\times\R^d)$ space when $g$ is restricted to the set $\mathcal{X}$. Here $f$ is the solution to the system
\begin{align} \label{eq:fix-point}
	\left\{ \begin{array}{ll}
		v \cdot \nabla_x f + E_g \cdot \nabla_v f = 0,  & \hbox{ in }\Omega\times \R^d,\\
		-\Delta \phi =  N - \rho_g,   & \hbox{ in } \Omega,\\ 
		\phi |_{\del \Omega} = 0 ,\\
		f |_{\Gamma_-} = \psi,\\
	\end{array}\right.
\end{align}
with $E_g = \nabla \phi$. Given the regularity of $g$, the characteristics (thus the map $\CalT$) are well-defined.

The proof is split into the following four steps. 

\smallskip
\noindent{\bf Step 1. $f$ is compactly supported in $\Sigma_\epsilon \times B(p_0, 2\epsilon)$.} 
This can be achieved by properly choosing sufficiently small $\epsilon$. 

We first solve the Vlasov equation, which is a pure transport equation, along the characteristics to get
\begin{align}\label{solution f}
   f(x, v) 
= \begin{cases}
     \psi(x_-(x, v), v_-(x, v)), & \text{if $(x, v)$ is on a characteristic curve initiating from $(x_-, v_-) \in \Gamma_-$}, \\
     0, & \text{otherwise}. 
   \end{cases} 
\end{align}
Next, for each $(x, v) \in \Supp(\psi)$, recall that by~\eqref{bound:t-plus}, the exit time of the corresponding characteristic satisfies the uniform bound
\begin{align} \label{bound:t-plus-main}
  t_+ (x, v) 
\leq 
  \frac{4 R}{|p_0|}, 
\qquad
  R = \text{diam}(\Omega),
\end{align}
provided that $|p_0| \geq 1 + 2\sqrt{2R \norm{E_g}_{L^\infty(\Omega)}}$. A sufficient condition is 
\begin{align} \label{bound:p-0-first}
   |p_0| 
\geq
  1 + 2\sqrt{2 R C\vpran{\norm{N}_{L^\infty(\Omega)} + \Eps^d  \norm{\psi}_{L^\infty_{x, v}}}}
\end{align}
due to Lemma~\ref{lemma:E} and $g\in \mathcal{X}$ ($\|g\|_{L^\infty_{x,v}}\leq \|\psi\|_{L^\infty_{x,v}}$). 
Applying the bound of $t_+$ to the integral equation~\eqref{eq:integ-charact} with $E_g$ as the given field, the difference between the X-ray ($x+vt$) and the characteristic ($X(t)$) satisfy
\begin{align}\label{EST:xt}
  |X(t) - x - v t|
\leq
  \frac{t_+^2}{2} \norm{E_g}_{L^\infty(\Omega)}
\leq
   \frac{8 C R^2}{|p_0|^2} \vpran{\norm{N}_{L^\infty(\Omega)} + \Eps^d \norm{\psi}_{L^\infty_{x, v}}}\leq \frac{8 C R^2}{|p_0|^2}(m_0+1),
\end{align}
where $C$ depends on $\Omega$ and $d$. For the velocity support of $f$, we use the integral form
\begin{align*} 
   V(t) = v + \int_0^t E_g(X(s)) \ds, 
\end{align*}
which leads to
\begin{align}\label{EST:vt}
   \abs{V(t) - v}  
\leq t_+ \norm{E_g}_{L^\infty(\Omega)}
\leq \frac{4CR}{|p_0|} \vpran{\norm{N}_{L^\infty(\Omega)} + \epsilon^d  \norm{\psi}_{L^\infty_{x, v}}}\leq \frac{4CR}{|p_0|} \vpran{m_0 + 1}. 
\end{align}
We then take 
\begin{align}\label{def:epsilon-0}
    \epsilon = \frac{c_0 }{|p_0|} \quad\hbox{ with }c_0 =  4CR \vpran{m_0 + 1}.
\end{align}
When $|p_0|$ is large enough such that ${2R\over |p_0|}\leq 1$, the estimates \eqref{EST:xt} and \eqref{EST:vt} satisfy
\begin{align*} 
   \abs{X(t) - x-vt}  \leq \epsilon\quad \hbox{ and }\quad \abs{V(t) - v}  \leq\epsilon 
\end{align*}
for all $(x,v)\in \Supp(\psi)$.
This implies the solution $f$ will be supported on $\Sigma\times B(p_0, 2\epsilon)$ based on \eqref{solution f}.

\smallskip
\noindent{\bf Step 2. $f$ is in  $W^{1, \infty}(\Omega \times \R^d)$.} In this step, we will show $f \in W^{1, \infty}(\Omega \times \R^d)$ and find its bounds. It is clear that $\norm{f}_{L^\infty_{x,v}} \leq \norm{\psi}_{L^\infty_{x, v}}$ due to \eqref{solution f}. Moreover, the differentiability of $f$ follows from that of $(x_-, v_-)$ and 
\begin{align*}
  \abs{\nabla_x f}
\leq
  \norm{\nabla_{x} \psi}_{L^\infty_{x, v}}
  \norm{\nabla_x x_-}_{L^\infty_{x, v}}
  + \norm{\nabla_{v} \psi}_{L^\infty_{x, v}}
     \norm{\nabla_x v_-}_{L^\infty_{x, v}}, 
\end{align*}
\begin{align*}
  \abs{\nabla_v f}
\leq
  \norm{\nabla_{x} \psi}_{L^\infty_{x, v}}
  \norm{\nabla_v x_-}_{L^\infty_{x, v}}
  + \norm{\nabla_{v} \psi}_{L^\infty_{x, v}}
     \norm{\nabla_v v_-}_{L^\infty_{x, v}}, 
\end{align*}
where bounds of the right-hand sides will follow from Lemma~\ref{lemma:Estimates}. 
To this end, consider the function
\begin{align*}
   h(x, v) = n_x \cdot \frac{v}{|v|}, 
\qquad
  (x, v) \in \del\Omega \times \R^d. 
\end{align*}
Since $\Omega$ is a smooth domain, the normal vector $n_x$ at $x\in\partial\Omega$ is Lipschitz continuous, that is, $\abs{n_x - n_{x_0}} \leq C(|x-x_0|)$ for some constant $C>0$. Then for each $(x, v) \in \Supp(\psi)$, 
\begin{align*}
  \abs{h(x, v) - h(x_0, p_0)}
\leq
  \abs{n_x - n_{x_0}} + \abs{\frac{v}{|v|} - \frac{p_0}{|p_0|}}
\leq
  C(|x-x_0|) + \frac{2|v-p_0|}{|p_0|}
\leq
  C \epsilon + \frac{2 }{|p_0| }\epsilon. 
\end{align*}
Therefore, for $|p_0|$ large enough, from \eqref{def:epsilon-0}, we obtain 
\begin{align} \label{bound:p-0-1}
  C \epsilon + \frac{2}{|p_0|}\epsilon \leq C{1\over |p_0|}
\leq
  \frac{1}{2} \delta_0. 
\end{align}
Together with \eqref{delta 0}, it yields that
\begin{align*}
   n_x \cdot \frac{v}{|v|} = h(x,v)\leq h(x_0,p_0) + {1\over 2}\delta_0 = - \frac{\delta_0}{2}   < 0
\qquad
  \text{for all $(x, v) \in \Supp(\psi)$}.  
\end{align*}
As a summary, we can choose the parameters in the bounds in Lemma~\ref{lemma:Estimates} as
\begin{align} \label{def:param-actual}
   \beta_0 = 4R, 
\qquad
  \delta = \frac{1}{2} \delta_0, 
\qquad
  M_0 = \norm{E_g}_{L^\infty(\Omega)},
\qquad
  M_1 = \norm{\nabla E_g}_{L^\infty(\Omega)}. 
\end{align}
Since $g \in \CalX$, by Lemma~\ref{lemma:E}, we have
\begin{align} \label{bound:M-0-1}
  M_0 
&\leq
  C\vpran{\norm{N}_{L^\infty(\Omega)} + \Eps^d  \norm{g}_{L^{\infty}(\Omega\times\R^d)}}
\leq
  C\vpran{\norm{N}_{L^\infty(\Omega)} + \Eps^d  \norm{\psi}_{L^{\infty}_{x,v}}},  \nn
\\
  M_1 
&\leq
  C\vpran{\norm{N}_{W^{1,\infty}(\Omega)} + \Eps^d  \norm{g}_{W^{1,\infty}(\Omega\times\R^d)}}
\leq 
  C \vpran{\|N\|_{W^{1,\infty}(\Omega)} + \epsilon^d (\|\psi\|_{L^\infty_{x,v}}+M_2)}. 
\end{align}
Now applying the explicit bounds \eqref{EST:diff x} - \eqref{EST:v diff x2} in Lemma~\ref{lemma:Estimates} and
using the bound of $\nabla_{x, v} \psi$ in~\eqref{bound:deriv-psi}, we have
\begin{align*}
 \norm{\nabla_x f}_{L^\infty_{x, v}}
+ \norm{\nabla_v f}_{L^\infty_{x, v}}
&\leq
   c_0' |p_0|  
  \vpran{\norm{\nabla_{x, v} x_-}_{L^\infty_{x, v}} 
             + \norm{\nabla_{x, v} v_-}_{L^\infty_{x, v}}}
\\
&\leq
  \CalH_1 
  + \CalH_2 \EE \, e^{\frac{C(1+M_1)}{|p_0|}}, 
\end{align*}
where we have 
\begin{align*}
  \CalH_1
&= C  |p_0| \vpran{1 + \frac{1}{\delta_0} + \frac{M_0}{\delta_0 |p_0|}}\leq  C|p_0|^2, 
\\
  \CalH_2
&= C  
   \vpran{1+\frac{1}{\delta_0|p_0|}+ {M_0 \over \delta_0|p_0|^2}}\leq C,
\end{align*}
by noting that $|p_0|<1$ and \eqref{bound:p-0-1}.
Inserting the bounds for $M_1$ in~\eqref{bound:M-0-1} and rearranging terms, we have
\begin{align*}
    \norm{\nabla_x f}_{L^\infty_{x, v}}
+ \norm{\nabla_v f}_{L^\infty_{x, v}}
&\leq
  \CalH_1
  + \CalH_3 
  \left(\|N\|_{W^{1,\infty}(\Omega)} + {c_0^d \over |p_0|^d} (\|\psi\|_{L^\infty_{x,v}}+M_2) \right) 
    e^{\frac{C c_0^d M_2 }{|p_0|^{d+1}} },
\end{align*}
 where
 \begin{align*}
   \CalH_3
 = C e^{C\left({1+ \|N\|_{W^{1,\infty}(\Omega)}\over |p_0|} + {c_0^d \|\psi\|_{L^\infty_{x,v}}\over |p_0|^{d+1}}\right)} 
    \CalH_2\leq  C. 
 \end{align*}
To choose the suitable $M_2$, we only need to set 
\begin{align*}  
  M_2 
\geq   \CalH_1
  + \CalH_3 
  \left(\|N\|_{W^{1,\infty}(\Omega)} + {c_0^d \over |p_0|^d}(\|\psi\|_{L^\infty_{x,v}}+M_2) \right)  
  e^{\frac{C c_0^d M_2 }{|p_0|^{d+1}} }.
\end{align*}
Since $d\geq 2$, we observe that the right-hand side is at most $O(|p_0|^2)$ when
\begin{align}\label{def:M-2}
    M_2= |p_0|^{2+\sigma}, \quad 0<\sigma <1
\end{align}
provided that $|p_0|$ is sufficiently large. With this choice of $M_2$, we conclude that $\norm{\nabla_{x, v} f}_{L^\infty_{x,v}}  \leq \GG$.
 
To summarize, by choosing $p_0$ as above, we obtain that
for $f = \CalT g \in \CalX$ for $g \in \CalX$.

\smallskip
\noindent{\bf Step 3. Lipschitz continuity of $\mathcal{T}$ in $L^\infty$.} 
We will show that $\CalT: \CalX \to \CalX$ is a contraction map given that $|p_0|$ is large enough. To this end, let $g_1, g_2 \in \CalX$ and $f_1 = \CalT g_1, f_2 = \CalT g_2$. Let $\tilde f = f_1 - f_2$. Then $\tilde f$ satisfies
\begin{align} \label{eq:tilde-f}
	\left\{ \begin{array}{ll}
		v \cdot \nabla_x \tilde f + E_{g_1} \cdot \nabla_v \tilde f
        = - \vpran{E_{g_1} - E_{g_2}} \cdot \nabla_v f_2,  & \hbox{ in } \Omega\times \R^d,\\
		\tilde f \big|_{\Gamma_-} = 0. 
	\end{array}\right.
\end{align}
The difference $E_{g_1} - E_{g_2}$ satisfies
\begin{align*}
   \norm{E_{g_1} - E_{g_2}}_{L^\infty(\Omega)}
=  \norm{\nabla_x\int_\Omega  G(x,y)(\rho_{g_1} - \rho_{g_2})(y)\,dy} 
\leq
   C \norm{\rho_{g_1} - \rho_{g_2}}_{L^\infty(\Omega)}, 
\end{align*}
where by the common support of $g_1$ and $g_2$ in $v$, we have
\begin{align*}
  \abs{\rho_{g_1} - \rho_{g_2}}
\leq 
  \int_{B(p_0, 2 \Eps)} \abs{g_1 - g_2} \dv
\leq
  c_d \Eps^d \norm{g_1 - g_2}_{L^\infty_{x,v}}.
\end{align*}
Hence, 
\begin{align*}
   \norm{E_{g_1} - E_{g_2}}_{L^\infty(\Omega)}
\leq
  \frac{C c_0^d }{|p_0|^d} \norm{g_1 - g_2}_{L^\infty_{x,v}}.
\end{align*}
Solving~\eqref{eq:tilde-f} along the characteristics gives
\begin{align*}
   \norm{f_1 - f_2}_{L^\infty_{x,v}}
&\leq
  \frac{C c_0^d}{|p_0|^d} \sup_{\Omega\times\R^d}t_+(x,v) \norm{\nabla_v f_2}_{L^\infty_{x,v}}\norm{g_1 - g_2}_{L^\infty_{x,v}}
\\
&\leq
   \frac{C c_0^{d}}{|p_0|^{d+1}} M_2 \norm{g_1 - g_2}_{L^\infty_{x,v}}, 
\end{align*}
where $t_+$ is the exit time associated with $E_{g_1}$. By the bound of $M_2$ in~\eqref{def:M-2}, it is clear that if we choose $|p_0|$ large enough such that
\begin{align} \label{bound:p-0-3}
   \lambda:=\frac{C c_0^{d+1}}{|p_0|^{d+1}} M_2
< 1,
\end{align}
then $\CalT$ is a contraction operator on $\mathcal{X}$ with respect to the $L^\infty$-norm.   

\smallskip
\noindent{\bf Step 4. Existence and uniqueness.} 
By iteration, for $g\in \mathcal{X}$, we let $f^0=g$ and define a sequence $f^{n+1}=\mathcal{T} f^{n}$ for $n\geq 0$. From step 1 and step 2, we have $\mathcal{T}f^n\in \mathcal{X}$ and from step 3, for $m>n\geq 1$, we can deduce
\begin{align*}
   \norm{ f^{m} - f^{n}}_{L^\infty_{x,v}}\leq {\lambda\over 1-\lambda} \norm{f^1 - f^0}_{L^\infty_{x,v}}, 
\end{align*}
which implies $\{f^n\}_{n\geq 1}$ is a Cauchy sequence. By completeness of $L^\infty$, there exists a limit $f\in L^\infty$ so that $f^n\rightarrow f$ in $L^\infty$ norm and $\|f\|_{L^\infty_{x,v}}\leq \|\psi\|_{L^\infty_{x,v}}$. Moreover, since $f^n\in \mathcal{X}$, for $(x,v),\,(y,w)\in\Omega\times\R^d$, we have
$$
    |f(x,v)-f(y,w)|=\lim_{n\rightarrow\infty} |f^n(x,v)-f^n(y,w)|\leq M_2 |(x,v)-(y,w)|,
$$
which implies $\|\nabla_{x,v}f\|_{L^\infty_{x,v}}\leq M_2$ and thus $f\in\mathcal{X}$.

Finally, the Lipschitz continuity of $\mathcal{T}$ yields
$$
    f = \lim_{n\rightarrow\infty} f^{n+1} = \lim_{n\rightarrow\infty} \mathcal{T}f^n = \mathcal{T}f.
$$
The uniqueness follows since such fixed point $f\in \mathcal{X}$ is unique.
\end{proof}

\section{Reconstruction of the doping profile}\label{sec:inverse problems} \label{sec:reconstruction}
With Theorem~\ref{thm:well-posedness}, we can now reconstruct the doping profile $N$. We denote 
    $$
        E_f = \nabla_x G \ast (N-\rho_f)=\nabla_x G \ast N-\nabla_x G \ast \rho_f=:\tilde E_N-\tilde{E}_\rho. 
    $$
Inspired by the work \cite{Novikov}, the following result states the asymptotic formula of the line integral of $\tilde{E}_N$, which can be recovered from the exit information on the boundary of the domain. 
\begin{thm}\label{thm:main-inverse}
    Let $\Omega$ be an open bounded and strictly convex domain with a $C^2$ boundary. Suppose $N \in \mathcal{M}$. For $(x_0,p_0)\in\Gamma_-$ with sufficiently large $|p_0|$, let $f$ be the solution to the system \eqref{eq:main} with boundary data $\psi\in \mathcal{S}_{x_0,p_0}(\Gamma_-)$. 
    Let $(X(t;x_0,p_0),V(t;x_0,p_0))$ be the corresponding characteristics (defined in \eqref{eq:ODE} with $E=E_f$) and let $x_+$, $v_+$ and $t_+$ be its forward exit point, velocity and time, respectively.
    
    Then $\tilde{E}_N$ satisfies
    \begin{align} \label{est:E-N-leading}
    \int_0^{|p_0| t_+^\ast} \tilde{E}_N(x_0 + e_{p_0} t) \dt
    = \lim_{|p_0| \to \infty} \vpran{v_+ - p_0}{|p_0|}, 
    \end{align}
    where $e_{p_0} = p_0/|p_0|$ and $t_+^*$ is the exit time of the line $X^*(t)=x_0+ p_0 t $ so that $(X^*(t^*_+), p_0 )\in \Gamma_+$.
\end{thm}
 
\begin{rmk}  
Since the line $X^*(t_+^*)=x_0+ p_0 t_+^* = x_0+e_{p_0}(|p_0|t_+^*)$, the upper limit $|p_0| t_+^\ast$ in \eqref{est:E-N-leading} is the total travel time of the line within $\Omega$ initiated from $x_0\in\partial\Omega$ with unit velocity $e_{p_0}$. Therefore, $|p_0| t_+^\ast$ depends only on the traveling distance and is independent of $|p_0|$. Equivalently, it denotes the total arc length of the line within $\Omega$.
\end{rmk}
 
 
\begin{proof} 
\smallskip
 
Let $(x_0,p_0)\in \Gamma_-$. To simplify the notation later on, let $t_+$ and $t_+^*$ be the forward exit times of trajectories $X(t)$ and $X^*(t)$ initiated at $(x_0,p_0)$, respectively.  
First we estimate the difference between $t_+$ and $t_+^\ast$. 
Recall that $\xi$ is the boundary defining function. We have
\begin{align} \label{eq:t-plus-and-ast}
      \xi(X(t_+;x_0,p_0)) = \xi\vpran{x_0 + p_0 t_+ + \int_0^{t_+} \int_0^s  E_f(X(\tau;x_0,p_0 )) \dtau \ds} = 0, 
\quad \text{and} \quad
      \xi\vpran{x_0 + p_0 t_+^\ast} = 0,
\end{align}
where $E_f = \tilde E_N - \tilde E_{\rho}$. 
Denote
\begin{align*}
	\tilde \eta := \int_0^{t_+} \int_0^s  E_f(X(\tau;x_0,p_0)) \dtau \ds. 
\end{align*}
By the bound of $E_f$ and $t_+$, we have
\begin{align*}
	\abs{\tilde \eta}
	\leq
	\frac{C}{|p_0|^2}
	\vpran{\norm{N}_{L^\infty(\Omega)} + \norm{\psi}_{L^\infty(\Omega)}}^3 =: \zeta_{p_0},
\end{align*}
where $C$ only depends on $\Omega$. It satisfies $|\tilde \eta| \ll 1$ if $|p_0|$ is large. 
Indeed, if we view $t_+$ and subsequently $\eta$ as given, then the equations in~\eqref{eq:t-plus-and-ast} can be generalized to the following expressions: 
\begin{align} \label{eq:t-plus-and-ast-1}
	\xi\vpran{x_0 + \tilde \eta + p_0 t_+^*(x_0+\eta,p_0)} = 0, 
	\quad \text{and} \quad
	\xi\vpran{x_0 + p_0 t_+^\ast} = 0,\quad |\eta| \ll 1. 
\end{align}
Extending the definition of the exit times to the points initiating outside of $\Omega$ for X-ray trajectories, we interpret the first equation in~\eqref{eq:t-plus-and-ast-1} as the $X$-ray $X^*$ starting from the point $(x_0+\eta, p_0)$ instead of the actual trajectory initiated at $(x_0, p_0)$. 
In particular, since $(x_0+\tilde{\eta},p_0)\in B_{2\zeta_{p_0}}(x_0,p_0)$, one has
$$
	t^*_+(x_0+\tilde{\eta},p_0) = t_+ (x_0, p_0)
$$
by observing that $X(t_+)=x_0 + \tilde \eta + p_0 t_+^*(x_0+\tilde\eta,p_0) \in\partial\Omega$. 
 Define a function $h:\R^n\times\R^n\times \R\rightarrow\R$ by $h(x,v,t):=\xi(x+vt)$. Then 
$$
    {\partial h\over \partial t}(x,v,t) = v\cdot \nabla \xi(x+vt).
$$
Due to the strict convexity of $\Omega$, there exists $\delta_1>0$ such that
$$
    e_{p_0}\cdot n_{x_+^*}>\delta_1>0,  
$$
and, moreover, for large enough $|p_0|$, we have
$$
    e_{p_0}\cdot n_{x+p_0t^*_+(x,p_0) } > {\delta_1\over 2}>0 \quad \hbox{ for } x\in B_{2\zeta_{p_0}}(x_0 ).  
$$	
Thus, by applying the implicit function theorem to~\eqref{eq:t-plus-and-ast-1}, we have for $x\in B_{\zeta_{p_0}}(x_0 )$,
$$
    \abs{\p_{x}t^*_+(x,p_0)} = \abs{{\p_{x}h(x,p_0,t^*_+(x,p_0))\over \p_t h(x,p_0,t^*_+(x,p_0))}}
    = \abs{ \frac{\nabla\xi(x+p_0t^*_+(x,p_0))}{p_0 \cdot \nabla\xi(x+p_0t^*_+(x,p_0))} } = {1\over|p_0|} {1\over e_{p_0}\cdot n_{x+p_0t^*_+(x,p_0) } }
    <
    \frac{2}{\delta_1 |p_0|}, 
$$
which leads to 
\begin{align}\label{bound:t-t-ast}
    \abs{t_+ - t_+^\ast}=|t^*_+(x_0+\tilde\eta,p_0)-t^*_+(x_0,p_0)|\leq |\p_{x}t^*_+(x,p_0)||\tilde\eta| = \frac{c_{N, \psi}}{\delta_1} \frac{1}{|p_0|^3}
    = \frac{o(1)}{|p_0|}, 
\end{align}
where $c_{N, \psi}$ depends on the $L^\infty$-bounds of $N, \psi$.

Now we are ready to show \eqref{est:E-N-leading} by comparing the two trajectories $X(t)$ and $X^*(t)$. We separate the two cases where $t_+ > t_+^\ast$ and $t_+ < t_+^\ast$. First, if $t_+ > t_+^\ast$, then 
\begin{align*}
   X(t) - x_0 - p_0 t = \int_0^t \int_0^s E_f (X(\tau)) d\tau ds 
\qquad
  t \in [0, t_+^\ast].
\end{align*}
Lemma~\ref{lemma:E} gives
\begin{align} \label{diff:X-L-1}
   \norm{X(t) - x_0 - p_0 t}_{L^\infty(\Omega)}
\leq
  \frac{C}{|p_0|^2} ,
\qquad
t \in [0, t_+^\ast],
\end{align}
where $C$ depends on $\Omega$, $d$, $N$, and $\psi$. 
Since $E_f= \tilde E_N - \tilde E_\rho$. we have
\begin{align*}
   v_+ - p_0 
&= \int_0^{t_+} E_f (X(\tau)) \dtau\\
&= \int_0^{t_+^*} \tilde{E}_N(X(\tau)) \dtau
   - \int_0^{t_+^*} \tilde{E}_\rho(X(\tau)) \dtau
   + \int_{t_+^\ast}^{t_+} E_f(X(\tau)) \dtau=:I_1+I_2+I_3. 
\end{align*}
We approximate the three terms on the right-hand side separately. The third term $I_3$ is of higher order since
\begin{align*}
  |I_3| = \abs{\int_{t_+^\ast}^{t_+} E_f(X(\tau)) \dtau}
\leq
  \norm{E_f}_{L^\infty(\Omega)} \abs{t_+ - t_+^\ast}
\leq
  C {o(1)\over |p_0|}
\end{align*}
by utilizing \eqref{bound:t-t-ast}.
To estimate the second term, we use the support of $f$ in $v$ to obtain that
\begin{align*}
   \norm{\tilde{E}_\rho}_{L^\infty(\Omega)}
\leq
  C\norm{\rho_f}_{L^\infty(\Omega)} =\norm{\int_{B(p_0,2\epsilon)}f\,dv}_{L^\infty(\Omega)} 
\leq
  \frac{C }{|p_0|^d},
\qquad
  d \geq 2,
\end{align*}
by recalling the definition $\epsilon$ in \eqref {def:epsilon-0}, where $C$ only depends on $\Omega$ and $\psi$.
Then the second term satisfies
\begin{align*}
  |I_2| = \abs{\int_0^{t_+^*} \tilde{E}_\rho(X(\tau)) \dtau}
\leq
  t_+^\ast \norm{\tilde{E}_\rho}_{L^\infty(\Omega)}
\leq
  \frac{C }{|p_0|^{d+1}}
\leq
  \frac{C }{|p_0|^{3}}, 
\qquad
  d \geq 2,
\end{align*}
where we used the fact that $t_+^*\leq {C\over |p_0|}$.
Finally, for the first term, we have
\begin{align*}
  I_1 = \int_0^{t_+^*} \tilde{E}_N(X(\tau)) \dtau
= \int_0^{t_+^*} \tilde{E}_N(x_0 + p_0 \tau) \dtau
   + \int_0^{t_+^*} \vpran{\tilde{E}_N(X(\tau)) - \tilde{E}_N(x_0 + p_0 t)} \dtau, 
\end{align*}
where by~\eqref{diff:X-L-1}, 
\begin{align*}
  \abs{\int_0^{t_+^*} \vpran{\tilde{E}_N(X(\tau)) - \tilde{E}_N(x_0 + p_0 \tau)} \dtau}
\leq
  \norm{\nabla \tilde{E}_N}_{L^\infty(\Omega)} 
  \int_0^{t_+^*} \abs{X(\tau) - (x_0 + p_0 \tau)} \dtau
\leq
  \frac{C }{|p_0|^3},
\end{align*}
where $C $ depends on $\Omega$, $N$, and $\psi$. Combining the estimates above, we obtain that if $t_+ > t_+^\ast$, then
\begin{align*}
  v_+ - p_0 
= \int_0^{t_+^*} \tilde{E}_N(x_0 + p_0 \tau) \dtau
   + {o(1)\over |p_0|}+\BigO\vpran{\frac{1}{|p_0|^3}}. 
\end{align*}
The case $t_+ < t_+^\ast$ can be dealt with similarly. In this case, 
\begin{align*}
   X(t) - x_0 - p_0 t = \int_0^t \int_0^s E_f (X(\tau)) d\tau ds, 
\qquad
  t \in [0, t_+],
\end{align*}
which gives
\begin{align} \label{diff:X-L-0}
   \norm{X(t) - x - p_0 t}_{L^\infty(\Omega)}
\leq
  {C \over |p_0|^2}, 
\qquad
t \in [0, t_+].
\end{align}
Then
\begin{align*}
   v_+ - p_0 
&= \int_0^{t_+} E_f (X(\tau)) \dtau
 = \int_0^{t_+} \tilde{E}_N(X(\tau)) \dtau
   - \int_0^{t_+} \tilde{E}_\rho(X(\tau)) \dtau
\\
& = \int_0^{t_+} \tilde{E}_N(x_0 + p_0 \tau) \dtau
      + \int_0^{t_+} \vpran{\tilde{E}_N(X(\tau)) - \tilde{E}_N(x_0 + p_0 \tau)} \dtau
      - \int_0^{t_+} \tilde{E}_\rho(X(\tau)) \dtau
\\
& = \int_0^{t_+^\ast} \tilde{E}_N(x_0 + p_0 \tau) \dtau
      - \int_{t_+}^{t_+^\ast} \tilde{E}_N(x_0 + p_0 \tau) \dtau
\\
& \quad   
      + \int_0^{t_+} \vpran{\tilde{E}_N(X(\tau)) - \tilde{E}_N(x_0 + p_0 \tau)} \dtau
      - \int_0^{t_+} \tilde{E}_\rho(X(\tau)) \dtau.
\end{align*}
The last three terms can be estimated similarly as for the case when $t_+ > t_+^\ast$. Overall, for both cases we have
\begin{align*}
    v_+ - p_0 
    = \int_0^{t_+^*} \tilde{E}_N(x_0 + p_0 \tau) \dtau
    + {o(1)\over |p_0|}+ \BigO\vpran{\frac{1}{|p_0|^3}}. 
\end{align*}
Thus, multiplying $|p_0|$ on both sides and taking the limit give
\begin{align*}
    \lim_{|p_0|\rightarrow \infty}|p_0|(v_+ - p_0) 
    = \lim_{|p_0|\rightarrow \infty}|p_0|\int_0^{t_+^*} \tilde{E}_N(x_0 + p_0 \tau) \dtau.
\end{align*}
By the change of variables $t=|p_0| \tau$, we have
\begin{align*}
    \lim_{|p_0|\rightarrow \infty}|p_0|(v_+ - p_0) 
    = \lim_{|p_0|\rightarrow \infty}\int_0^{|p_0| t_+^\ast} \tilde{E}_N(x_0 + e_{p_0} t) \dt= \int_0^{|p_0| t_+^\ast} \tilde{E}_N(x_0 + e_{p_0} t) \dt,
\end{align*}
which gives the leading-order approximation of the integrals of $\tilde{E}_N$ over $\Omega$ as in~\eqref{est:E-N-leading}. Note that $|p_0| t^\ast_+$ is independent of $|p_0|$.
 \end{proof}

\begin{rmk}
In Theorem~\ref{thm:main-inverse}, $v_+=v_+(x_0,p_0)$ is the exit velocity that corresponds to $(x_0, p_0)\in \Gamma_-$. 
To extract $v_+$ from the Albedo operator, one can choose the incoming data $\psi$ such that $0 \leq \psi \leq 1$ and satisfies $\psi(x,v) = 1$, $(x,v)\in\Gamma_-$ if and only if $(x,v) = (x_0,p_0)$. Since the Vlasov equation does not change the amplitude of the density along each trajectory, the solution can be formally written as $f(x,v) = \psi(x_-(x,v),v_-(x,v))=\delta_{x_0,p_0}(X(-t_-(x,v)),V(-t_-(x,v)))$ with Dirac delta function $\delta_{x_0,p_0}$, 
where $x_-(x,v), \, v_-(x,v),\,t_-(x,v)$ are the backward exit point, velocity and time of $(x,v)\in \overline\Omega\times\R^d$, respectively. Hence $(x_+(x_0,p_0), v_+(x_0,p_0))$ can be determined through $f(x_+(x_0,p_0), v_+(x_0,p_0)) = 1$. We refer to a similar observation in \cite{McDowallCPDE} where the scattering relation, the set of all incoming and outgoing points and directions, is determined from knowledge of the Albedo operator on the boundary. 
\end{rmk}

\subsection{Proof of Theorem~\ref{thm:uniqueness} (Uniqueness result)}
\begin{proof}[Proof of Theorem~\ref{thm:uniqueness}]
    From Theorem~\ref{thm:main-inverse}, we can deduce that for any $(x_0,p_0)\in \Gamma_-$,
    \begin{align*}
        \int_0^{|p_0| t_+^\ast}  (\tilde{E}_{N_1}-\tilde{E}_{N_2})(x_0+e_{p_0}t)\,dt =0 
    \end{align*}
    with total travel time $|p_0| t_+^\ast$ in $\Omega$.
    By extending $\tilde{E}_{N_1}-\tilde{E}_{N_2}$ by zero outside $\Omega$, this implies that the X-ray transform of $\tilde{E}_{N_1}-\tilde{E}_{N_2}$ vanishes, namely,
    $$
        R(\tilde{E}_{N_1}-\tilde{E}_{N_2})(x_0,e_{p_0}):=\int_\R (\tilde{E}_{N_1}-\tilde{E}_{N_2}) (x_0+e_{p_0}t)\,dt =0,
    $$
    where the X-ray transform of a function $h$ is defined by $Rh(x,\theta)=\int_\R h(x+\theta s)\,ds$.
    Now since $\tilde{E}_{N_1}-\tilde{E}_{N_2}\in L^1(\R^d)$ has compact support,
    the injectivity of the X-ray transform (see, for example, [\cite{Joonas}, Corollary~1.4]) yields $\tilde{E}_{N_1}(x)=\tilde{E}_{N_2}(x)$ in $\Omega$. This implies 
    $$
       \text{div}\tilde{E}_{N_1}(x)=\text{div} \tilde{E}_{N_2}(x).
    $$
    Moreover, since 
    $\tilde{E}_{N_j} = \nabla_x G\ast N_j$, which is the gradient of the solution $u_j$ to $\Delta u=N_j$ with $u|_{\partial\Omega}=0$, we derive that 
\begin{equation*}
       0 = \text{div}\tilde{E}_{N_1}-\text{div} \tilde{E}_{N_2} = \Delta u_1 - \Delta u_2 = N_1-N_2.   \qedhere
\end{equation*}
\end{proof}

\begin{rmk} \label{rmk:generalization}
Note that our analysis uniquely recovers the electric field $E_N$ generated by the doping profile, which is given by 
\begin{align*}
   E_N = \nabla_x G \ast N,
\end{align*}
from which we have recovered N. In the case where $N$ is known and $G$ is unknown, our analysis can be applied to recover the interaction potential $G$, which is a problem of interest for material sciences. 
\end{rmk}


\bibliography{ref}
\bibliographystyle{abbrv}

\end{document}